\documentclass{article}
\pdfoutput=1
\usepackage{spconf}
\usepackage[named]{algo}
\usepackage{graphics}
\usepackage{graphicx}
\usepackage{stfloats}
\usepackage{amsmath}
\usepackage{amsfonts}
\usepackage{subfigure}
\usepackage{float}
\usepackage{url}

\newcommand{\Section}{\section}
\newcommand{\SubSection}{\subsection}

\newcommand{\bPsi}{\mbox{\boldmath   $\Psi$}}

\def\urltilda{\kern -.15em\lower .7ex\hbox{\~{}}\kern .04em}
\def\urldot{\kern -.10em.\kern -.10em}
\def\urlhttp{http\kern -.10em\lower -.1ex\hbox{:}\kern -.12em\lower 0ex\hbox{/}\kern -.18em\lower 0ex\hbox{/}}

\ninept

\begin{document}

\title{A Fast Algorithm for the Constrained Formulation
of Compressive Image Reconstruction and Other Linear
Inverse Problems}

\name{Manya V. Afonso, Jos\'{e} M. Bioucas-Dias, and M\'{a}rio A. T. Figueiredo}
\address{ Instituto de Telecomunica\c{c}\~{o}es, \\
Instituto Superior T\'{e}cnico, Technical University of Lisbon, {\bf Portugal}\\
Email: $\{$mafonso, jose.bioucas, mario.figueiredo$\}$@lx.it.pt\thanks{M. Afonso
is supported by a EU Marie-Curie Fellowship (EST-SIGNAL program: {\scriptsize\tt est-signal.i3s.unice.fr});
contract MEST-CT-2005-021175.}}
\maketitle

\begin{abstract}
Ill-posed linear inverse problems (ILIP), such as restoration and reconstruction, are a core topic  of signal/image processing. A standard approach to deal with ILIP uses a constrained optimization problem, where a regularization function is minimized under the constraint that the  solution explains the observations sufficiently well. The regularizer and constraint are usually convex; however, several particular features of these problems (huge dimensionality, non-smoothness) preclude the use of off-the-shelf optimization tools and have stimulated much research. In this paper, we propose a new efficient  algorithm to handle one class of constrained problems (known as basis pursuit denoising) tailored to image recovery applications. The proposed algorithm, which belongs to the category of augmented Lagrangian methods, can be used to deal with a variety of imaging ILIP, including deconvolution and reconstruction from compressive observations (such as MRI). Experiments testify for the effectiveness of the proposed method.
\end{abstract}


\Section{Introduction}
\SubSection{Problem Formulation}
\label{sec:intro}
Linear inverse problems constitute one of the central themes of
signal/image processing. In this class of problems, a noisy
indirect observation ${\bf y}$, of an original signal ${\bf x}$,
is modeled as
\begin{equation}
{\bf y} = {\bf B}{\bf x} + {\bf n},\label{observation_model}
\end{equation}
where ${\bf B}$  is the matrix representation of the direct
operator and ${\bf n}$ is noise. In the sequel, we denote by $n$
the dimension of ${\bf x}$, thus ${\bf x}\in\mathbb{R}^n$,
while  ${\bf y}\in\mathbb{R}^m$.
In the classical problem of image deblurring/deconvolution, ${\bf B}$
is the matrix representation of a convolution operator. In other
reconstruction problems, ${\bf B}$ represents some linear direct operator,
such as of tomographic projections (Radon transform) or a partially
observed (e.g., Fourier) transform (as in compressive MRI \cite{Lustig}).

Usually, the problem of estimating ${\bf x}$ from ${\bf y}$
is ill-posed ({\it e.g.}, if $m < n$),
thus requiring some sort of regularization.
In the presence of noise, a natural criterion to infer
${\bf x}$ from  ${\bf y}$ has the form \cite{ChambolleLions,Malgouyres}
\begin{equation}
\min_{{ \bf x}} \phi(\bf x)\hspace{0.5cm} \mbox{subject to} \hspace{0.25cm}
 \|{\bf B}{\bf x}-{\bf y}\|_{2} \leq \varepsilon,
\label{genconstrained}
\end{equation}
where $\phi: \mathbb{R}^n \rightarrow \bar{\mathbb{R}}$ is the the regularizer
and $\varepsilon \geq 0$ a parameter which depends on the noise variance.
In the case where $\phi({\bf x}) = \|{\bf x}\|_1$, the above problem is usually known
as {\it basis pursuit denoising} (BPD) \cite{ChenDonohoSaunders}. The
{\it basis pursuit} (BP) problem is the particular case of (\ref{genconstrained})
for $\varepsilon = 0$. In recent years, an explosion of interest in problems
of the form  (\ref{genconstrained}) was sparked by the emergence of
{\it compressive sensing} (CS) \cite{Candes}, \cite{donoho}. The theory of
CS  provides conditions (on matrix ${\bf B}$ and the
degree of sparseness of the original ${\bf x}$) under which a solution
of (\ref{genconstrained}) is an optimal (in some sense) approximation
to the ``true" ${\bf x}$.

In most image recovery and CS problems, the regularizer $\phi$ is convex but
non-smooth; typical examples are the {\it total variation} (TV)
\cite{Candes}, \cite{art:Rdin:O:F:Physica:92} and  $\ell_1$ norms.
Problem (\ref{genconstrained}) is thus convex, but the very
high dimension (usually $\geq 10^4$) of ${\bf x}$ and ${\bf y}$
precludes the direct application of off-the-shelf optimization
algorithms. This difficulty is further amplified by the fact that matrix ${\bf B}$
only ``exists" as an operator; {\it i.e.}, there are efficient algorithms
to compute products of ${\bf B}$ (or ${\bf B}^T$) by some vector (image),
but it is highly impractical to extract and manipulate individual blocks,
rows, or columns of this matrix.

\subsection{Previous Work}
Most state-of-the-art methods for dealing with linear inverse
problems, under convex, non-smooth regularizers (namely, TV and $\ell_1$),
consider, rather than (\ref{genconstrained}), the unconstrained problem
\begin{equation}
\min_{{\bf x}} \frac{1}{2}\| {\bf B\, x} -
{\bf y}\|_2^2 + \tau \, \phi ({\bf x}),\label{unconstrained}
\end{equation}
where $\tau \in \mathbb{R}_+$ is the so-called regularization parameter.
Of course, problems (\ref{genconstrained}) and (\ref{unconstrained})
are equivalent, in the following sense: for any $\varepsilon$ such that
problem (\ref{genconstrained}) is feasible, a solution of (\ref{genconstrained})
is either the null vector, or else is a solution of
(\ref{unconstrained}), for some $\tau$  \cite{GPSR}.

The currently fastest (publicly available) algorithms for
solving (\ref{unconstrained}), include:
gradient projection for sparse reconstruction (GPSR)
 \cite{GPSR};  fast iterative shrinkage/thresolding algorithm (FISTA) \cite{FISTA};
 two-step IST (TwIST) \cite{TwIST}; and sparse reconstruction by
 separable approximation (SpaRSA) \cite{SpaRSA_SP}.
These methods were shown to
be considerably faster than earlier methods, including $l1\verb=_=ls$
\cite{KimKLBG07} and the codes in the $\ell_1$-{\it magic}
({\tt http://www.l1-magic.org}) and  the {\it SparseLab} ({\tt http://sparselab.stanford.edu})
toolboxes. Very recently, we have introduced a new algorithm,
called SALSA ({\it split augmented Lagrangian shrinkage algorithm});
 experiments on a set of standard image recovery problems show
that SALSA is faster than GPSR, TwIST, FISTA, and SpaRSA \cite{FigueiredoDiasAfonso}.

Although it is usually easier/simpler to solve an
unconstrained problem than a constrained one,
formulation (\ref{genconstrained}) has an important
advantage: parameter $\varepsilon$ has a clear
meaning (it is proportional to the noise variance) and is much
easier to set than parameter $\tau$ in (\ref{unconstrained}).
Of course, one may solve (\ref{genconstrained}) by using
one of the algorithms mentioned in the previous paragraph to solve
(\ref{unconstrained}) and searching for the ``correct" value of
$\tau$ that makes (\ref{unconstrained}) equivalent to
(\ref{genconstrained}). Clearly, this is not efficient,
as it involves solving many instances of (\ref{unconstrained}).
Obtaining fast algorithms for solving (\ref{genconstrained}) is thus
an important research front.

There are few efficient algorithms to solve (\ref{genconstrained})
in an image recovery context: ${\bf x}$ and ${\bf y}$  of dimension
$\geq 10^4$ (often $\geq 10^6$), ${\bf B}$ representing an operator,
and $\phi$ a convex, non-smooth function. A notable exception is the recent
SPGL1 \cite{SPGL1}, which (as its name implies) is specifically
designed for $\ell_1$ regularization ($\phi({\bf x}) = \|{\bf x}\|_1$).
Other methods for solving problems with the form (\ref{genconstrained}),
for $\phi$ equal to the $\ell_1$ or TV norms, are available in the
$\ell_1$-magic package; however, as shown in \cite{SPGL1}, those methods are
 quite inefficient for large problems. General purpose  methods, such as
those in the SeDuMi package ({\tt http://sedumi.ie.lehigh.edu}), are simply not applicable
when ${\bf B}$ is not an actual matrix, but an operator.

The Bregman iterative algorithm (BIA) was recently proposed to solve
(\ref{genconstrained}) with $\varepsilon = 0$,  but is not
directly applicable when $\varepsilon > 0$ \cite{YinOsherGoldfarbDarbon}.
To deal with the case of $\varepsilon > 0$, it was suggested
that the BIA for $\varepsilon = 0$ is used and stopped when
$\|{\bf B}{\bf x}-{\bf y}\|_{2}\leq \varepsilon$ \cite{CaiOsherZhen},
\cite{YinOsherGoldfarbDarbon}. Clearly, that approach is not guaranteed to find a good
solution, since it depends strongly on the initialization; {\it e.g.},
if the algorithm starts at a feasible point, it will immediately stop,
although the point may be far from a minimizer of $\phi$.

\SubSection{Proposed Approach}

In this paper, we introduce an algorithm for solving
optimization problems of the form (\ref{genconstrained}). The
basic ingredients are the following: the original
constrained problem (\ref{genconstrained}) is transformed
into an unconstrained one by using an indicator function of
the feasible set; the resulting unconstrained problem is
transformed into a different constrained problem, by the application
of a variable splitting operation; finally, the
obtained constrained problem is attacked with an {\it
augmented Lagrangian} (AL) scheme \cite{NocedalWright},
which is a variant of the {\it alternating direction method of multipliers} (ADMM)
\cite{EcksteinBertsekas}. Since (as SALSA), the proposed
method uses variable splitting and AL optimization,
we call it C-SALSA (for {\it constrained}-SALSA).

The resulting algorithm is more general than SPGL1, in
the sense that it can be used with any convex regularizer
$\phi$ for which the corresponding Moreau proximity
operator \cite{CombettesSIAM}, defined as
\begin{equation}
\bPsi_{\tau\phi}({\bf y}) = \arg\min_{{\bf x}} \frac{1}{2}\|{\bf x}-{\bf y}\|_2^2 + \tau\phi({\bf x}),
\label{MPM}
\end{equation}
has closed form or can be efficiently computed.
Below, we will show examples of C-SALSA where ${\bf x}$
is an image, $\phi$ is the TV norm \cite{art:Rdin:O:F:Physica:92},
and $\bPsi_{\tau\phi}$ is computed using Chambolle's
algorithm \cite{Chambolle}. Another classical choice is
$\phi({\bf x})=\|{\bf x}\|_1$, which leads to
$\bPsi_{\tau\phi}({\bf y}) = \mbox{soft}({\bf y}, \tau)$, where
$\mbox{soft}(\cdot, \tau)$ denotes the component-wise
application of the soft-threshold function $y \mapsto \mbox{sign}(y)\max\{|y|-\tau,0\}$.

C-SALSA is experimentally shown to efficiently solve
image recovery problems of the form (\ref{genconstrained}),
such as MRI reconstruction from CS-type partial Fourier
observations using TV regularization.
Moreover, C-SALSA is also shown to be faster than
SPGL1 in wavelet-based image deconvolution
problems under  $\ell_1$ regularization.

The paper is organized as follows. Section \ref{sec:tools}
briefly reviews variable splitting and ADMM. Section \ref{sec:salsa}
contains the derivation leading to C-SALSA.
Section \ref{sec:experiments} reports experimental results, and
Section \ref{sec:conclusions} ends the paper with a few remarks and
pointers to future work.

\Section{Variable Splitting and ADMM}
\label{sec:tools}

Consider an unconstrained optimization problem
\begin{equation}
\min_{{\bf u}\in \mathbb{R}^n} f_1({\bf u}) +
f_2\left({\bf G}{\bf u}\right),\label{unconstrained_basic}
\end{equation}
where ${\bf G} \in \mathbb{R}^{d\times n}$.
Variable splitting (VS) is a  simple procedure that consists in
creating new variables, say ${\bf v}$ and ${\bf w}$,
to serve as the argument of each of the terms, $f_1$ and $f_2$,
under the constraints that ${\bf w} = {\bf u}$ and $ {\bf v} = {\bf G}{\bf u}$, that is,
\begin{equation}
 {\displaystyle \min_{{\bf u},{\bf w}\in \mathbb{R}^n ,\, {\bf v}\in\mathbb{R}^d }}\;\;  f_1({\bf w}) + f_2({\bf v}),
 \hspace{0.6cm} \mbox{subject to} \hspace{0.3cm} \begin{array}[t]{l} {\bf w} = {\bf u}\\ {\bf v} = {\bf G}{\bf u}.
 \end{array} \label{constrained_basic}
\end{equation}
Problem (\ref{constrained_basic}) is clearly equivalent to
the unconstrained problem (\ref{unconstrained_basic}).
The rationale behind VS  is that it may
be easier to solve the constrained problem (\ref{constrained_basic})
than to solve its unconstrained counterpart (\ref{unconstrained_basic}).
It is important to stress that the VS in (\ref{constrained_basic}) is not
the one commonly used, where only variable ${\bf v}$ is created; however,
as shown below, the proposed VS will lead to a very effective algorithm.

Other variants of VS were recently used in several image processing problems:
in \cite{Wang}, it was used to obtain a fast TV-based algorithm;
in \cite{BioucasFigueiredo2008}, it was used to handle problems with compound
regularization. VS also underlies the recent split-Bregman methods \cite{GoldsteinOsher},
but there the splitting is different and with a different goal.

Using an augmented Lagrangian (AL) approach to handle  problem (\ref{constrained_basic})
leads to the following algorithm, also known as the {\it method of multipliers} (MM)
\cite{Hestenes}, \cite{Powell}  (see also \cite{FigueiredoDiasAfonso}, for details):
\begin{eqnarray}
\left({\bf u}_{k+1} , {\bf v}_{k+1}, {\bf w}_{k+1} \right) & \in & \arg\min_{{\bf u},{\bf v},{\bf w}}\;\left\{ f_{1}({\bf w})
 + f_{2}({\bf v})  + \rule[-0.2cm]{0cm}{0.4cm} \right.\nonumber \\
\lefteqn{\hspace{-1cm} \left.
\rule[-0.2cm]{0cm}{0.4cm} \frac{\mu_{1}}{2} \|{\bf G}{\bf u} - {\bf v} - {\bf b}_k\|_2^2
  +\frac{\mu_{2}}{2} \|{\bf u} - {\bf w} - {\bf c}_k\|_2^2\right\}}\label{mixed}\\
{\bf b}_{k+1} & = & {\bf b}_{k} + {\bf G}{\bf u}_{k+1} - {\bf v}_{k+1}\\
{\bf c}_{k+1} & = & {\bf c}_{k} + {\bf u}_{k+1} - {\bf w}_{k+1}.
\end{eqnarray}

Problem (\ref{mixed}) is not trivial since it
involves non-separable quadratic as well as non-smooth terms.
Replacing (\ref{mixed}) by the alternating minimization with
respect to each vector leads to a variant of the so-called
{\it alternating direction method of multipliers} (ADMM)
\cite{EcksteinBertsekas}:
\begin{algorithm}{ADMM}{
\label{alg:salsa3}}
Set $k=0$, choose $\mu_1,\mu_2 > 0$, ${\bf v}_0$, ${\bf w}_0$, ${\bf b}_0$, and ${\bf c}_0.$\\
\qrepeat\\
$  {\bf u}_{k+1}  \leftarrow  {\displaystyle \arg\min_{{\bf u}}}
\frac{\mu_{1}}{\mu_2} \|{\bf G}{\bf u} \! - \!{\bf v}_k \! -\! {\bf b}_k\|_2^2 + \|{\bf u} \!- \!{\bf w}_k \! - {\bf c}_k\|_2^2$\\
  $  {\bf v}_{k+1}  \leftarrow  {\displaystyle \arg\min_{{\bf v}}} f_{2}({\bf v})
 + \frac{\mu_{1}}{2} \|{\bf G}{\bf u}_{k+1} - {\bf v} - {\bf b}_k\|_2^2$\\
   $  {\bf w}_{k+1}  \leftarrow  {\displaystyle \arg\min_{{\bf w}}} f_{1}({\bf w})
 + \frac{\mu_{2}}{2} \|{\bf w} - {\bf u}_{k+1} - {\bf c}_k\|_2^2$\\
     ${\bf b}_{k+1} \leftarrow {\bf b}_{k} + {\bf Gu}_{k+1} - {\bf v}_{k+1}$\\
     ${\bf c}_{k+1} \leftarrow {\bf c}_{k} + {\bf u}_{k+1} - {\bf w}_{k+1}$\\
     $k \leftarrow k+1$
\quntil stopping criterion is satisfied.
\end{algorithm}

The proof of convergence in \cite{EcksteinBertsekas} applies to a
different variant of ADMM, which results from a different splitting.
However, it is possible to show that this version can
still be written as a standard ADMM and satisfies the conditions of the
convergence theorem \cite{FigueiredoDias2009}.

\Section{Proposed Method}
\label{sec:salsa}
\vspace{-0.1cm}
\SubSection{Reformulation of the Problem}
\vspace{-0.1cm}
The feasible set in problem (\ref{genconstrained}) is the ellipsoid
\begin{equation}
E(\varepsilon,{\bf B},{\bf y}) =
\{{\bf x}\in \mathbb{R}^n:\| {\bf B\, x}- {\bf y}\|_2 \leq \varepsilon\},
\end{equation}
possible infinite in some directions. Problem (\ref{genconstrained}) can be written
as an unconstrained (discontinuous) problem,
\begin{equation}
\min_{{\bf x}} \; \phi({\bf x}) + \iota_{E(\varepsilon,{\bf I},0)} ({\bf B x - y}) ,
\label{unconstrained_reformulation}
\end{equation}
where $\iota_S:\mathbb{R}^m\rightarrow \bar{\mathbb{R}}$
denotes the indicator function of set $S\subset \mathbb{R}^m$,
\begin{equation}
\iota_S({\bf s})=\left\{
\begin{array}{ll}
0, & \text{if } {\bf s} \in S\\
+\infty, & \text{if } {\bf s} \notin S.
\end{array} \right.\label{indicator_ellipsoid}
\end{equation}
Notice that $E(\varepsilon,{\bf I},0)$ is simply an $\varepsilon$-radius Euclidean ball
centered at the origin of  $\mathbb{R}^m$.

Since problem (\ref{unconstrained_reformulation}) clearly has the
form (\ref{unconstrained_basic}), its VS-based constrained optimization reformulation
is
\begin{equation}
 {\displaystyle \min_{{\bf u},{\bf w}\in \mathbb{R}^n, {\bf v}\in \mathbb{R}^m}} \; \phi({\bf w}) +
 \iota_{E(\varepsilon,{\bf I},0)} ({\bf v}),
 \hspace{0.4cm} \mbox{s. t.} \hspace{0.15cm} \begin{array}[t]{l} {\bf u} = {\bf w}\\
{\bf v} = {\bf B u - y}.\end{array}
\label{constrained_image_linear}
\end{equation}

\vspace{-0.2cm}
\SubSection{Application of ADMM}
\vspace{-0.2cm}
Performing the adequate translations (which are clear from
comparing (\ref{constrained_basic}) with (\ref{constrained_image_linear})),
the ADMM becomes the proposed C-SALSA.

{\footnotesize
\begin{algorithm}{C-SALSA}{
\label{alg:salsa4}}
Set $k=0$, choose $\mu_1,\mu_2 > 0$, ${\bf v}_0$, ${\bf w}_0$, ${\bf b}_0$, and ${\bf c}_0.$\\
\qrepeat\\
   ${\bf u}' \leftarrow {\bf w}_{k} + {\bf c}_k$ \\
   ${\bf u}'' \leftarrow  {\bf y} + {\bf v}_{k} + {\bf b}_k$ \\
   ${\bf u}_{k+1} \leftarrow {\displaystyle \arg\min_{{\bf u}}} \; \frac{\mu_1}{\mu_2} \|{\bf B\, u} - {\bf u}''\|_2^2
   + \|{\bf u}-{\bf u}' \|_2^2 $\\
   ${\bf v}' \leftarrow {\bf Bu}_{k+1} - {\bf y} - {\bf b}_k$\\
   ${\bf v}_{k+1} \leftarrow {\displaystyle  \arg\min_{{\bf v}}} \; \iota_{E(\varepsilon,{\bf I},0)} ({\bf v}) + \frac{\mu_1}{2}\|{\bf v}-{\bf v}'\|_2^2$\\
   ${\bf w}' \leftarrow {\bf u}_{k+1} + {\bf c}_k$\\
   ${\bf w}_{k+1} \leftarrow {\displaystyle  \arg\min_{{\bf w}}} \; \phi({\bf w}) + \frac{\mu_2}{2}\|{\bf w}-{\bf w}'\|_2^2$\\
  ${\bf b}_{k+1} \leftarrow {\bf b}_{k} + {\bf Bu}_{k+1} - {\bf y} - {\bf v}_{k+1}$\\
     ${\bf c}_{k+1} \leftarrow {\bf c}_{k} + {\bf u}_{k+1} - {\bf w}_{k+1}$\\
     $  k \leftarrow k+1$
\quntil stopping criterion is satisfied.
\end{algorithm}}

A key feature of C-SALSA is that the cost of each iteration is
$O(n\log n)$, as confirmed by the following observations.
Lines 3, 4, 8, 11, and 12 simply involve adding vectors or scalars,
thus have $O(n)$ or $O(1)$ cost. Line 5 consists in
minimizing a strictly convex quadratic function, leading  (with $\alpha = \mu_1/\mu_2$) to
\begin{equation}
{\bf u}_{k+1} = \left(\alpha \, {\bf B}^T {\bf B} + {\bf I}\right)^{-1}
\left( \alpha \, {\bf B}^T {\bf u}'' +  {\bf u}' \right).\label{eq:Wiener}
\end{equation}
As will be shown in Subsection \ref{sec:computingxk}, in several cases of interest,
this matrix inversion has $O(n\log n)$ cost. Lines 6 and 10 involve
matrix-vector products which, by the same reason, have $O(n\log n)$ cost.
Line 7 corresponds to the
orthogonal projection of ${\bf v}'$ onto the
$\varepsilon$-radius $\ell_2$ ball $E(\varepsilon,{\bf I},0)$, which is an $O(n)$ operation:
\begin{equation}
{\bf v}_{k+1} = {\cal P}_{E(\varepsilon,{\bf I},0)}({\bf v}') =
 \left\{ \begin{array}{ll}
\varepsilon \, {\bf v}' /  \|{\bf v}'\|_2, & \text{if }\;  \|{\bf x}\|_2 > \varepsilon,\\
{\bf v}', & \text{if }\;  \|{\bf v}\|_2 \leq \varepsilon.
\end{array} \right.\label{projection_ellipsoid}
\end{equation}
Finally, line 9 is simply ${\bf w}_{k+1} = \bPsi_{\phi/\mu_2}({\bf w}')$
(see (\ref{MPM})). If $\phi({\bf x}) = \|{\bf x}\|_1$, the cost of
$\bPsi$ is $O(n)$. If $\phi$ is the TV norm, we use Chambolle's
algorithm, which (although iterative) also has $O(n)$ cost  \cite{Chambolle}.

\SubSection{Implementing (\ref{eq:Wiener})}
\label{sec:computingxk}
We will now show how (\ref{eq:Wiener}) can be implemented with
$O(n\log n)$ cost in several cases of interest.
If ${\bf B}$ represents a convolution,  it is  factorized as
${\bf B = U}^H {\bf D U}$, where ${\bf U}$ is the unitary matrix
(${\bf U}^H = {\bf U}^{-1} $) representing the discrete Fourier transform (DFT)
and ${\bf D}$ is a diagonal matrix. Thus,
\begin{equation}
(\alpha {\bf B}^T{\bf B} + {\bf I})^{-1} =  {\bf U}^H\left(\alpha |{\bf D}|^2 + {\bf I}\right)^{-1}{\bf U},\label{eq:Wiener2}
\end{equation}
where $|{\bf D}|^2$ is the  matrix with squared absolute values
of the entries of ${\bf D}$. Since $\alpha |{\bf D}|^2 + {\bf I}$ is diagonal,
its inversion costs $O(n)$. Products by ${\bf U}$  and ${\bf U}^H$ have
$O(n\log n)$ cost, using the FFT algorithm.

In frame-based regularization, the unknown image is represented
on a frame ({\it e.g.}, of wavelets or curvelets)
and then the coefficients of this representation are estimated
from the observed data, under some regularizer. A constrained
formulation of this approach still has the form (\ref{genconstrained})
but with different meanings for ${\bf x}$ and ${\bf B}$: vector
${\bf x}$ now contains the frame coefficients
of the unknown image ${\bf Wx}$ (the columns of ${\bf W}$
contain the elements of the adopted frame) and ${\bf B} = {\bf AW}$
is now the product of an observation  matrix ${\bf A}$ by the
frame synthesis matrix ${\bf W}$ \cite{SpaRSA_SP}. The only impact of
this change on C-SALSA is in computing (\ref{eq:Wiener}), since
${\bf AW}$ is not diagonalizable by the DFT.
This difficulty may be sidestepped under the assumption that
${\bf W}$ contains a 1-tight (Parseval) frame
({\it i.e.}, ${\bf W\,W}^H = {\bf I}$) and that ${\bf A = U}^H {\bf D U}$,
with ${\bf D}$ diagonal ({\it e.g.}, a convolution). Using the matrix inversion lemma:
\begin{equation}
\!\!\left(\alpha {\bf W}^H\!{\bf A}^H\!{\bf AW}\! + {\bf I}\right)^{-1} \!\!\! = \,
 {\bf I} - \! {\bf W}^H\!\overbrace{{\bf A\!}^H \!\left({\bf AA\!}^H + {\bf I}/{\alpha}\right)^{-1}\!\!\!\!{\bf A}}^{\bf F}{\bf W}.
\label{eq:filter_3}
\end{equation}
Since ${\bf A = U}^H {\bf D U}$, we have ${\bf F} = {\bf U}^H {\bf D^*} \left(|{\bf D}|^2  +  {\bf I}/\alpha\right)^{-1} {\bf D U}$,
the computation of which has $O(n\log n)$ cost, using  the FFT to compute the
products by ${\bf U}$ and ${\bf U}^H$.
The  cost of (\ref{eq:filter_3}) will thus be either $O(n\log n)$
or the cost of the products by ${\bf W}^H$ and ${\bf W}$. For most tight
frames used in image processing, there are fast $O(n\log n)$ algorithms
to compute these products \cite{Mallat}.

Finally, we considered the case of partial Fourier observations,
which is used to model  MRI  acquisition and has been the focus of
recent interest due to its connection
to compressed sensing \cite{Candes}, \cite{Lustig}.
In this case, ${\bf B}={\bf M}{\bf U}$, where ${\bf M}$ is an $m\times n$ binary
matrix ($m<n$) formed by a subset of rows of the identity, and ${\bf U}$ was defined above.
Due to its particular structure, matrix ${\bf M}$ satisfies
${\bf M}{\bf M}^T = {\bf I}$; this fact together with the matrix inversion lemma
leads to
\begin{eqnarray}
\left(\alpha{\bf B}^T{\bf B} + {\bf I}\right)^{-1} =
{\bf I} - \alpha/(1+\alpha)\, {\bf U}^H {\bf M}^T{\bf M}{\bf U},
\label{eq:mask_matrix}
\end{eqnarray}
where ${\bf M}^T{\bf M}$ is equal to an identity with some
zeros in the diagonal. Consequently,  the cost of (\ref{eq:mask_matrix})
is also $O(n\log n)$.

\Section{Experiments}
\label{sec:experiments}
\vspace{-0.2cm}
All experiments were performed using MATLAB on a Windows XP laptop with a $2$ GHz processor and $512$ MB of RAM.

\begin{figure}[ht]
\centering
\includegraphics[width=0.25\textwidth]{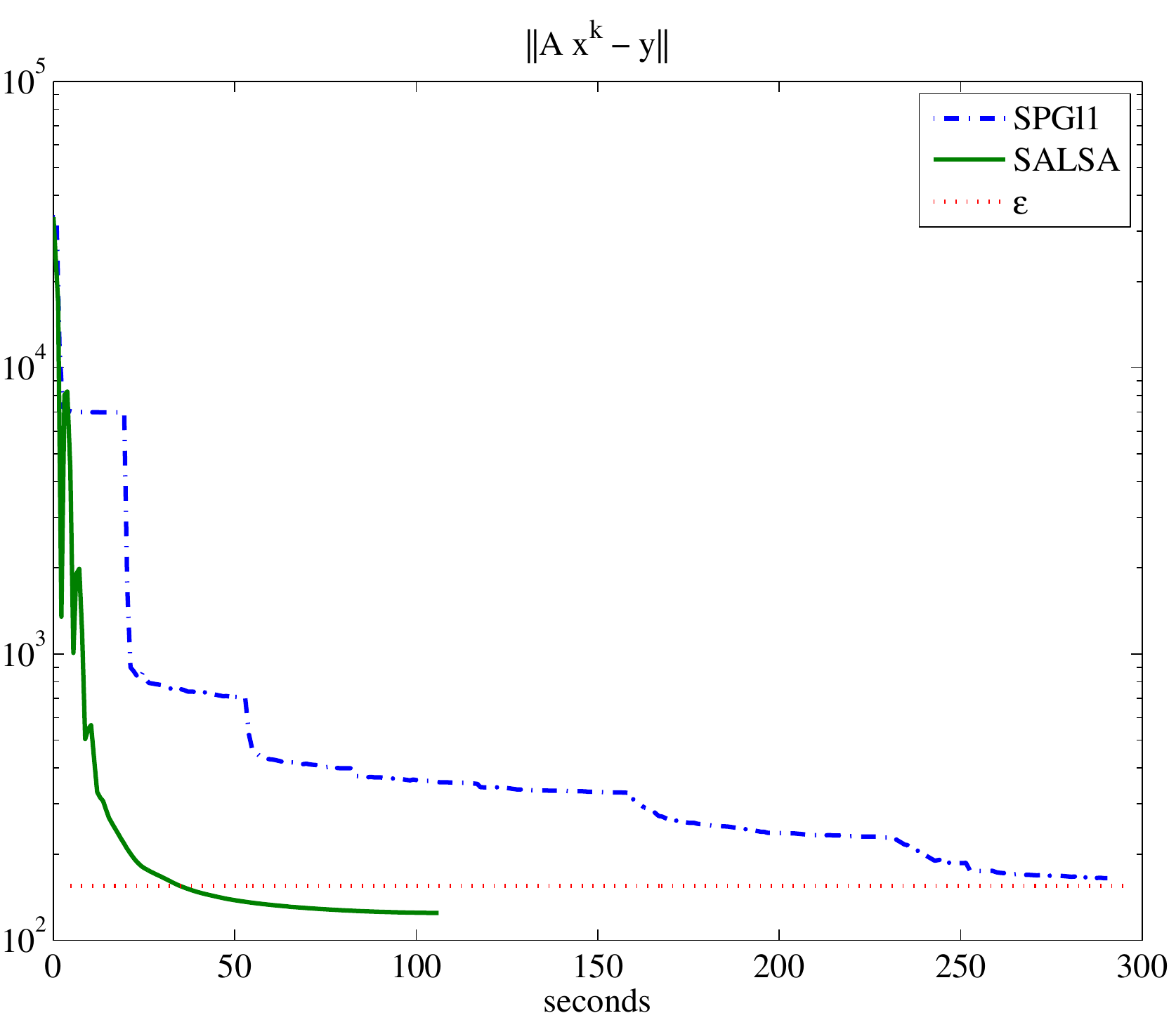}
\label{fig:evolutioncriterion}
\caption{\footnotesize Image deblurring with wavelets: Evolution of the quadratic constraint $\|{\bf AWx}-{\bf y}\|_2$ over time for $9\times 9$ uniform
blur, $\sigma=0.56$}
\end{figure}

We consider five standard image deconvolution benchmark problems \cite{FigueiredoNowak2003},
summarized in Table~\ref{decon_problems}, all on the well-known Cameraman image.
We solve problem (\ref{genconstrained}), with $\phi({\bf x})=\|{\bf x}\|_1$
 (thus $\bPsi$ is a soft threshold) and ${\bf B=AW}$, where ${\bf W}$ is a redundant 4-level Haar wavelet frame and ${\bf A}$ is the blur operator.
We set $\mu_1 = \mu_2$ and hand-tuned its value for fastest convergence. We compare C-SALSA
with SPGL1 as follows. First, we run SPGL1 and then C-SALSA (from the same initialization),
stopping when the constraint in (\ref{genconstrained}) is satisfied and the MSE of
the estimate is below that obtained by SPGL1. Table~\ref{tab:deconresultsl1redundant} reports
the number of iterations and CPU times taken in each of the experiments.
Figure~\ref{fig:evolutioncriterion} plots the evolution of quadratic constraint
$\|{\bf AWx}_k - {\bf y}\|_2$, in experiment $1$.

\begin{table}[h]
\centering
\caption{Details of the image deblurring experiments.}\label{decon_problems}
\begin{tabular}{|c|l|l|} \hline
\footnotesize  Experiment & \footnotesize blur kernel  \rule[-0.1cm]{0cm}{0.4cm}  & \footnotesize $\sigma^2$ \\ \hline
\footnotesize 1  & \footnotesize $9\times 9$ uniform & \footnotesize $0.56^2$ \\
\footnotesize 2A & \footnotesize Gaussian & \footnotesize 2\\
\footnotesize 2B & \footnotesize Gaussian & \footnotesize 8\\
\footnotesize 3A  & \footnotesize $h_{ij} = 1/(1 + i^2 + j^2)$ & \footnotesize 2 \\
\footnotesize 3B  & \footnotesize $h_{ij} = 1/(1 + i^2 + j^2)$ & \footnotesize 8 \\
\hline
\end{tabular}
\end{table}

\begin{table}[h]
\centering \caption{Image deblurring using wavelets - Computation speed}
\label{tab:deconresultsl1redundant}
\begin{tabular}{|l|c|c| c|c|}
\hline
\footnotesize Experiment & \multicolumn{2}{|c|}{\footnotesize Iterations} & \multicolumn{2}{|c|}{\footnotesize CPU time (seconds)}\\
\hline
 & \footnotesize SPGL1 &\footnotesize  C-SALSA & \footnotesize SPGL1 &\footnotesize  C-SALSA\\
\hline
\footnotesize 1 &\footnotesize  400 &\footnotesize  136 & \footnotesize 553.188 & \footnotesize 118.953\\
\footnotesize 2A &\footnotesize  200 &\footnotesize  152 & \footnotesize 258.406  & \footnotesize 130.203\\
\footnotesize 2B &\footnotesize  150 &\footnotesize  120 & \footnotesize 190.688 & \footnotesize 115.375\\
\footnotesize 3A &\footnotesize  250 &\footnotesize  57 & \footnotesize 303.688 & \footnotesize 48.5\\
\footnotesize 3B &\footnotesize  150 &\footnotesize  46 & \footnotesize 188.516 & \footnotesize 40.5156\\
\hline
\end{tabular}
\end{table}

In the MRI reconstruction experiment, ${\bf M}$ models 22 radial observations of the DFT and $\phi$ is the TV norm \cite{Candes}.
Since SPGL1 can be used only for $\phi({\bf x})=\|{\bf x}\|_1$, we compare
C-SALSA with the code available in $\ell_1$-magic. Table~\ref{tab:MRIresults}
compares the 2 algorithms, in terms of computation time and the final MSE obtained.

\begin{table}[h]
\centering \caption{MRI reconstruction - Comparison of computation speed}
\label{tab:MRIresults}
\begin{tabular}{|l|c| c|c|}
\hline
\footnotesize Algorithm  & \footnotesize CPU time (seconds) & \footnotesize MSE\\
\hline
\footnotesize  $\ell_1$-magic  &\footnotesize  710.997 & \footnotesize 0.000117224\\
\footnotesize C-SALSA  & \footnotesize 18.6875 &\footnotesize  6.79023e-007\\
\hline
\end{tabular}
\end{table}


\Section{Conclusions}
\label{sec:conclusions}
\vspace{-0.1cm}
We have proposed a fast algorithm for solving constrained convex
optimization problems usually known as {\it basis pursuit denoising}.
Our algorithm is based on variable splitting and exploits
augmented Lagrangian tools. Preliminary experiments with $\ell_1$ and TV
regularization show that the new algorithm
outperforms existing methods in terms of computation time, by a considerable
factor. Ongoing work includes a more thorough experimental evaluation
of C-SALSA.

\end{document}